\documentclass[12pt]{amsart}
\usepackage{amsfonts}

\theoremstyle{definition}

\theoremstyle{remark}

\newcommand{\const}{\mathop{\rm const}\limits}

\newcommand{\supp}{\mathop{\rm supp}\limits}

\newcommand{\mes}{\mathop{\rm mes}\limits}

\newcommand{\diag}{\mathop{\rm diag}\limits}

\let \vs = \vspace

\begin{document}

\begin{center}

{\bf EXACT NORM ESTIMATES FOR MULTIVARIATE DILATION  \\

\vspace{3mm}

  OPERATORS BETWEEN TWO BILATERAL  \\

\vspace{3mm}

  WEIGHT GRAND LEBESGUE SPACES}\par

\vspace{4mm}

{\bf E. Ostrovsky}\\

e-mail: eugostrovsky@list.ru \\

\vspace{3mm}

{\bf L. Sirota}\\

e-mail: sirota3@bezeqint.net \\

\vspace{4mm}

Department of Mathematics and Statistics, Bar-Ilan University, \\
 59200, Ramat Gan, Israel. \\

\vspace{5mm}

 Abstract. \\

\end{center}

\vspace{4mm}

{\it We give in this short paper  a sharp estimate for the norm
 of a  multivariate dilation  operator generated by
 multi-matrix (tensor) linear argument transformation (dilation operator)
 between two different weight Lebesgue-Riesz and Grand Lebesgue Spaces  (GLS). \par
 We consider also some examples and study the compactness  of these operators.}

\vspace{4mm}


\vspace{4mm}


\vspace{3mm}

\section{Introduction. Notations. Statement of problem.}

\vspace{3mm}

{\bf A. Matrix dilation in Lebesgue-Riesz spaces. } \par

\vs{3mm}

Let here $  X = R^d  $ equipped with Borelian sigma-field and with  Lebesgue
measure $  dx  $ and let $ f: R^d \to R  $ be some (measurable) function  belonging to the space $  L_p(R^d), \ p \ge 1. $ The
norm of a function $ f $ this space will be denoted as ordinary

$$
|f|L_p = |f|_p \stackrel{def}{=} \left[ \int_X |f(x)|^p \ dx \right]^{1/p}.
$$

 Let also $  A $ be non degenerate: $ \det(A) \ne 0  $ linear (homogeneous) map (matrix) acting from $  R^d $ to $  R^d. $ \par

 Define a {\it matrix  dilation}, or {\it compression}, in the terminology of an article \cite{Shimogaki1},  operator of a form

 $$
 V_A[f] = f(Ax). \eqno(1.1)
 $$
 Obviously,

$$
| V_A[f]|_p^p = \int_{R^d} |f(A x)|^p \ dx = \int_{R^d}  |\det(A)|^{-1} \  |f(y)|^p \   dy =
$$

$$
|\det(A)|^{-1} \ |f|_p^p,
$$
or equally

$$
|V_A[f]|_p = |\det(A)|^{-1/p}  \ |f|_p. \eqno(1.2)
$$

 In particular, if  $ d = 1 $ and following $  V_A[f](x) = f(A\cdot x), \ A = \const \ne 0; $ then $  V_A[\cdot] $
is the classical dilation operator, and we have

$$
|V_A[f]|_p = |A|^{-1/p}  \ |f|_p. \eqno(1.2a)
$$

 The last equalities (1.2) (and (1.2a)) may be rewritten as follows. Note that the fundamental function $  \phi(L_p(R^d), \delta) $
 of the $ L_p(R^d) $ space has a form

$$
  \phi(L_p(R^d), \delta)  = \delta^{1/p}, \delta \in (0,\infty).
$$

  At the  same result is true for arbitrary Lebesgue-Riesz space $  L_p, \ 1 \le p < \infty $ constructed over any
measure space equipped with diffuse measure, see \cite{Bennett1}, chapters 1,2. \par

 Recall that for arbitrary rearrangement invariant space $ (S, ||\cdot||S)  $ (over $ R^d $ or arbitrary  another measure space
 with diffuse measure) the fundamental function $ \phi(S, \delta)  $ is defined as follows

$$
\phi(S, \delta)  := ||I_D(\cdot)||S,
$$
where $  I_D (\cdot) $ denotes an indicator function of a measurable set $ D $ and $ \mes(D) = \delta, \ \delta \in [0, \infty]. $\par
 So,

$$
|V_A[f]|L_p =  \phi(L_p,|\det A|^{-1})  \cdot |f|L_p \eqno(1.3)
$$
or equally

$$
||V_A[f]||(L_p \to L_p) =  \phi(L_p,|\det A|^{-1}). \eqno(1.3a)
$$

\vspace{3mm}

{\bf B. Matrix dilation in weighted Lebesgue-Riesz spaces. } \par

\vspace{3mm}

Let $ |x|, \ x \in R^d $ be arbitrary (complete) norm in the whole space $  R^d,  $ for instance, the classical Euclidean, and
let $ ||A|| $ be  correspondent matrix norm:

$$
||A|| \stackrel{def}{=} \sup_{ 0 \ne x \in R^d} \left[ \frac{|Ax|}{|x|}  \right],
$$
so that $  |Az| \le ||A|| \cdot |z|,  z \in R^d.  $\par

Let also $ \alpha = \const > 0 $ and $ \mu_{\alpha}(\cdot)  $ be {\it weight} measure defined on the Borelian sets on $ R^d: $

$$
\mu_{\alpha}(D) \stackrel{def}{=} \int_D |x|^{\alpha} dx.
$$
 The correspondent for the space $ L_{p,\alpha} $ norm for the function $  f: R^d \to R $ will be denoted by $ |f|_{p,\alpha}: $

$$
|f|^p_{p,\alpha} := \int_{R^d} |f(x)|^p \ |x|^{\alpha} \ dx.
$$

We deduce

$$
|  V_A [f] |_{p,\alpha}^p = \int_{R^d} |f(A x)|^p \ |x|^{\alpha} \ dx = | \det A|^{-1} \int_{R^d} |f(y)|^p \ |A^{-1}y|^{\alpha} \ dy \le
$$

$$
 | \det A|^{-1} \ ||A^{-1}||^{\alpha} \int_{R^d} |f(y)|^p \ |y|^{\alpha} \ dy = | \det A|^{-1} \ ||A^{-1}||^{\alpha} \ |f|^p_{p,\alpha};
$$

$$
|  V_A [f] |_{p,\alpha} \le  | \det A|^{-1/p} \ ||A||^{-\alpha/p} \ |f|_{p,\alpha},
$$
or equally

$$
||V_A||(L_{p,\alpha} \to L_{p,\alpha}) \le | \det A|^{-1/p} \ ||A||^{-\alpha/p}. \eqno(1.4)
$$

 Equivalent form:

$$
||V_A||(L_{p,\alpha} \to L_{p,\alpha}) \le \phi\left(L_{p,\alpha}, |\det(A)|^{-1} \cdot ||A||^{-\alpha} \right). \eqno(1.5)
$$

  If for example the matrix $  A $ is diagonal: $ A = \diag(\lambda_1, \lambda_2, \ldots, \lambda_d), $ and following
 $  |\det A| = \prod_k |\lambda_k|,  $ then the relation (1.4) can be transformed as an equality:

$$
||V_A||(L_{p,\alpha} \to L_{p,\alpha}) = | \det A|^{-(1 + \alpha)/p}, \ \det A \ne 0. \eqno(1.6)
$$

 In particular, if $ \lambda_k = \lambda = \const \ne 0, \ k = 1,2, \ldots,d  $ then

$$
||V_A||(L_{p,\alpha} \to L_{p,\alpha}) = | \lambda|^{-d(1 + \alpha)/p}. \eqno(1.7)
$$

\vspace{3mm}

{\bf  Our purpose in this short article is investigation of these operators in Grand  Lebesgue Spaces (GLS),
as well as  in  its multivariate  anisotropic version AGLS. \par

We set ourselves a goal to derive the sharp estimates for the norms of these estimates, similar to the estimates (1.2), (1.7). } \par

\vspace{3mm}

  Analogous statement of this problem for another Banach spaces, namely: Lorentz, Orlicz, Marcinkiewicz  etc.
 closed relatively as a rule with computation of the so-called Boyd's and Shimogaki etc.
 indices and as a consequence with application in turn in the theory of operators interpolation, theory of
 Fourier series  in the one-dimensional case $ d = 1 $  see, e.g. in monograph and articles \cite{Bennett1}, chapter 7;
  \cite{Boyd1}, \cite{Boyd2}, \cite{Montgomery-Smith1}, \cite{Ostrovsky109}, \cite{Cornelia1}, \cite{Shimogaki1} etc. \par

 The immediate predecessor of offered article is the authors preprint \cite{Ostrovsky109}, where are obtained an upper bounds
for the norms of dilation operators acting in these spaces in the one - dimensional case.\par

 Another operators acting in GLS spaces: Hardy, Riesz, Fourier, maximal, potential, composition etc.  are investigated, e.g.  in
 \cite{Liflyand1}, \cite{Ostrovsky100}, \cite{Ostrovsky101}, \cite{Ostrovsky3}, \cite{Ostrovsky4}. \\

\bigskip

\section{Grand Lebesgue Spaces (GLS). }

\vspace{3mm}

We recall first of all here  for reader convenience  some definitions and facts from
the theory of GLS spaces.\par

 Recently, see
\cite{Fiorenza1}, \cite{Fiorenza2},\cite{Ivaniec1}, \cite{Ivaniec2}, \cite{Jawerth1},
\cite{Karadzov1}, \cite{Kozatchenko1}, \cite{Liflyand1}, \cite{Ostrovsky1}, \cite{Ostrovsky2} etc.
 appear the so-called Grand Lebesgue Spaces GLS

 $$
 G(\psi) = G = G(\psi ; a;b);  \ a;b = \const; \ a \ge 1, \ b \le \infty
 $$
spaces consisting on all the measurable functions $ f : X \to R  $ with finite norms

$$
||f||G(\psi) \stackrel{def}{=} \sup_{p \in (a;b)} \left[\frac{|f|_p}{\psi(p)} \right]. \eqno(2.1)
$$

 Here $ \psi = \psi(p), \ p \in (a,b) $ is some continuous positive on the {\it open} interval $ (a;b) $ function such
that

$$
\inf_{p \in(a;b)} \psi(p) > 0.
$$

 We define formally  $ \psi(p) = +\infty, \ p \notin [a,b]. $\par

We will denote
$$
\supp(\psi) \stackrel{def}{=} (a;b).
$$

The set of all such a functions with support $ \supp(\psi) = (A,B) $ will be denoted by  $  \Psi(A,B). $  \par

These spaces are non-trivial, are rearrangement invariant; and are used, for example, in
the theory of Probability, theory of Partial Differential Equations,
 Functional Analysis, theory of Fourier series,
 Martingales, Mathematical Statistics, theory of Approximation  etc. \par

They does not coincide in general case with  classical rearrangement invariant spaces: Lorentz, Orlicz,
Marcinkiewicz etc., see \cite{Kozatchenko1}, \cite{Ostrovsky1}, \cite{Liflyand1}, \cite{Ostrovsky7}. \par

 Notice that the classical Lebesgue - Riesz spaces $ L_p $  are extremal case of Grand Lebesgue Spaces;
the exponential Orlicz spaces are the particular cases of Grand Lebesgue Spaces, see
 \cite{Ostrovsky2},  \cite{Ostrovsky100}.  \par

  Let a function $  f:  X \to R  $ be such that

 $$
 \exists (a,b): \ 1 \le a < b \le \infty \ \Rightarrow  \forall p \in (a,b) \ |f|_p < \infty.
 $$
Then the function $  \psi = \psi(p) $ may be {\it naturally} defined by the following way:

$$
\psi_f(p) := |f|_p, \ p \in (a,b).
$$

\vspace{3mm}

 Recall now that the fundamental  function  $ \phi(G\tau, \ \delta), \ 0 \le \delta \le \mu(X)  $
 for the Grand Lebesgue Space   $  G \tau  $ may be calculated by the formula

$$
\phi(G\tau, \ \delta) = \sup_{p \in (a,b)} \left[ \frac{\delta^{1/p}}{\tau(p)}  \right].
$$

 This notion play a very important role in the theory of operators, Fourier analysis etc., see \cite{Bennett1}.
The detail investigation of the fundamental  function  for GLS is done in \cite{Liflyand1},  \cite{Ostrovsky2}. \par

\bigskip

\section{Main result: the norm of matrix dilation operator on GLS. }

\vspace{3mm}

 Let the measurable function $ f: R^d \to R  $ be such that there exists a function $ \psi(\cdot) $ from the set
 $ \Psi(a,b), 1 \le a < b \le \infty  $ for which $ f \in G\Psi(a,b).  $   For instance, the function $ \psi(\cdot)   $
 can be picked as a natural function for the function $  f: \ \psi(p) := |f|_p, $ if of course the last function is finite
 for all the values $  p  $ from some non-trivial interval $  (a,b). $ \par

 Suppose the function $ \psi(\cdot) $ has a form

$$
\psi(p) = \frac{\nu(p)}{\zeta(p)},  \ p \in (a,b), \eqno(3.0)
$$
where  both the functions $ \nu(\cdot), \ \zeta(\cdot) $ belong also to at the same set $ G\Psi(a,b).  $ \par

\vs{3mm}

{\bf Theorem 3.1.} Assume (in the notations of the first two sections) $ \det(A) \ne 0. $ Then

$$
||V_A[f]||G\nu \le \phi(G\zeta, |\det(A)|^{-1}) \ ||f||G\psi, \eqno(3.1)
$$
and the last estimate is in  general case non-improvable. \par

\vs{3mm}

{\bf Proof.} We can and will suppose without loss of generality $ ||f||G\psi = 1.  $ Then

$$
|f|_p \le \psi(p), \ p \in (a,b).
$$

 We use the inequality (1.2):

$$
|V_A[f]|_p = |\det(A)|^{-1/p}  \ |f|_p \le |\det(A)|^{-1/p}  \ \psi(p),
$$
or equally

$$
\frac{|V_A[f]|_p}{\nu(p)} \le \frac{|\det(A)|^{-1/p}}{\zeta(p)}. \eqno(3.2)
$$
 Let us take the supremum over $ p \in (A,B) $ from both the sides of the last inequality (3.2), taking into account
the direct definitions of norm and fundamental function for Grand Lebesgue Spaces:

$$
||V_A[f]||G\nu \le \phi\left(G\zeta, |\det(A)|^{-1} \right) = \phi\left(G\zeta, |\det(A)|^{-1} \right) \ ||f||G\psi. \eqno(3.3)
$$

Inequality (3.3) becomes an equality if for example $ \psi(p) = |f|_p,  $ i.e. the function $ \psi(p) $ is the natural function
for the source function $  f(\cdot). $\par
 This completes the proof of theorem 3.1. \par





\vs{3mm}

 \section{ Weight estimates for GLS.   }

\vspace{3mm}

  The notion of Grand  Lebesgue Spaces may be easy generalized on the arbitrary measure space with sigma-finite measure, see
\cite{Fiorenza1}-\cite{Ostrovsky2} etc. Denote for instance

$$
||f||G\psi_{\alpha} \stackrel{def}{=} \sup_{ p \in (a,b)} \left[ \frac{|f|_{p,\alpha}}{\psi(p)}  \right], \eqno(4.0)
$$
where as before $ \alpha = \const > 0, \ \psi(\cdot) \in \Psi(a,b), \ 1 \le a < b \le \infty.  $ \par

\vs{3mm}

Assume again in the notations of the first two sections $ x \in R^d, \ \det(A) \ne 0. $
Let also $ \psi(\cdot), \ \zeta(\cdot)  $ be two functions  from the set $  \Psi(a,b).  $
Introduce the new function $ \nu = \nu(p) $ by an equality

$$
\nu(p) := \psi(p) \cdot \zeta(p); \eqno(4.1)
$$
obviously, the function  $ \nu = \nu(p)  $ belongs also to the set $ \Psi(a,b). $

\vs{4mm}

{\bf Theorem 4.1.}

$$
||V_A[f]||G\nu_{\alpha}  \le \phi(G\zeta_{\alpha}, |\det(A)|^{-1} \ ||A||^{-\alpha} ) \ ||f||G\psi_{\alpha}, \eqno(4.2)
$$
and the last estimate is also in  general case non-improvable, for example  when  the matrix $  A  $ is diagonal and
the function $ \psi(\cdot)  $ is natural for the function $ f: \ \psi(p) = |f|_p < \infty, \ p \in (a,b).  $ \par

\vs{3mm}

{\bf Proof } is alike to one in theorem 3.1.  Indeed, let $ ||f||G\psi_{\alpha} = 1. $
We the inequality (1.4):

$$
|  V_A [f] |_{p,\alpha} \le  | \det A|^{-1/p} \ ||A||^{-\alpha/p} \ |f|_{p,\alpha} \le  | \det A|^{-1/p} \ ||A||^{-\alpha/p} \ \psi(p),
$$
or equally

$$
\frac{|  V_A [f] |_{p,\alpha} }{\nu(p)} \le \frac{ |\det A|^{-1/p} \ ||A||^{-\alpha/p}}{\zeta(p)}. \eqno(4.3)
$$
It remains to take the supremum from both the sides of the last inequality (4.3) over $ p: \ p \in (a,b): $

$$
||V_A[f]||G\nu_{\alpha}  \le \phi(G\zeta_{\alpha}, |\det(A)|^{-1} \ ||A||^{-\alpha} )  =
\phi(G\zeta_{\alpha}, |\det(A)|^{-1} \ ||A||^{-\alpha} ) \ ||f||G\psi_{\alpha},
$$
Q.E.D.\par

\vs{3mm}

 \section{ Main result:  dilation operators in mixed (anisotropic) Lebesgue spaces.}

\vspace{5mm}

  We recall here the definition of the so-called anisotropic (mixed in Bochner's sense)
Lebesgue (Lebesgue-Riesz) spaces; see the source work  \cite{Benedek1}. More detail information about this
spaces see in the classical books  of Besov O.V., Il’in V.P., Nikol’skii S.M.
\cite{Besov1}, chapter 16,17; Leoni G. \cite{Leoni1}, chapter 11; using for us theory of
operators interpolation in this spaces see in \cite{Besov1}, chapter 17,18. \par

  Let $ (X_j,A_j,\mu_j), \ j = 1,2,\ldots, l $ be measurable spaces with sigma-finite non - trivial measures $ \mu_j; $
in the considered  in this report case $ X_j = R^{m_j}. $  \par

 Set

$$
X = R^d  = \otimes_{j=1}^l X_j,
$$
evidently $  d = \sum_j m_j.  $ \par

 Let also

 $$
 p = \vec{p} = (p_1, p_2, . . . , p_l) \eqno(5.1)
 $$
  be $ l- $ dimensional numerical vector such that $ 1 \le p_j \le \infty.$ \par

 Recall that the anisotropic Lebesgue space $ L_{ \vec{p}} $ consists on all the  total measurable
real valued function  $ f = f(x_1,x_2,\ldots, x_l) = f( \vec{x} ) $

$$
f:  \otimes_{j=1}^l X_j \to R
$$

with finite norm $ |f|_{ \vec{p} } \stackrel{def}{=} $

$$
\left( \int_{X_l} \mu_l(dx_l) \left( \int_{X_{l-1}} \mu_{l-1}(dx_{l-1}) \ldots \left( \int_{X_1}
 |f(\vec{x})|^{p_1} \mu_1(dx_1) \right)^{p_2/p_1 }  \ \right)^{p_3/p_2} \ldots   \right)^{1/p_l}.
$$

 Note that in general case $ |f|_{p_1,p_2} \ne |f|_{p_2,p_1}, $
but $ |f|_{p,p} = |f|_p. $ \par

 Observe also that if $ f(x_1, x_2) = g_1(x_1) \cdot g_2(x_2) $ (condition of factorization), then
$ |f|_{p_1,p_2} = |g_1|_{p_1} \cdot |g_2|_{p_2}, $ (formula of factorization). \par

\vspace{3mm}

{\bf Definition 5.1.} Let $  D  $ be Borelian subset of the whole space $  R^d $ and

$$
 q = \vec{q} = (q_1, q_2, . . . , q_d) \eqno(5.2)
$$
 be $ d \ - $ dimensional numerical vector such that $ 1 \le q_j < \infty.$ \par
We define as before as a capacity of a {\it fundamental function} the following expression

$$
\phi_{\vec{q}}(D) \stackrel{def}{=} || I_D(\cdot) ||_{\vec{q}}. \eqno(5.3)
$$

\vs{3mm}

{\bf Remark 5.1.} \
 Note that in general case $ d \ge 2 $ the value $ \phi_{\vec{q}}(D) $ does not dependent only on the volume, i.e.
 on the measure of the set $  D. $\par

{\bf Remark 5.2.} \ This notion of fundamental function, or in other words, {\it  non-linear volume,  }
 may be easy generalized on the arbitrary multivariate, for instance, over the Euclidean
space $ R^d, $ Banach functional space $ (L, ||\cdot||L)  $ as follows

$$
\phi_L(D) \stackrel{def}{=}  || I_D(\cdot) ||L, \eqno(5.4)
$$
if there exists.

\vs{3mm}

 This function $  D \to \phi_L(D)  $  is obviously non-negative and sub-additive:

$$
\phi_L(\cup_{k=1}^n D_k) \le \sum_{k=1}^n \phi_L(D_k), \ 1 \le n \le \infty.
$$
 If the space $ L $ coincides with the classical $  L_1 $  ones, then the function $ D \to \phi_L(D) $ is ordinary sigma additive
Lebesgue measure. In the case when $  L = L_{p,p,p, \ldots,p}  $ we return to the fundamental function for $ L_p(R^d) $ space. \par

\vs{3mm}

{\bf Remark 5.3.} \ Assume that the set $  D  $ is direct (Cartesian) product of the (measurable) sets $  F_j: $

$$
D = \otimes_{j=1}^g F_j, \ F_j \subset R^{m_j}.
$$
 Since  the indicator function $ I_D $ is factorable, we deduce by means of the formula of factorization

$$
\phi_{\vec{q}}(D) = \prod_{j=1}^g \phi_{ \vec{q_j} } (F_j).  \eqno(5.5)
$$
 As a consequence: in this case the value  $ \phi_{\vec{q}}(D) $ dependent  only on the "individual" volumes
 $  \{ \mu_j(F_j)  \}. $ \par

\vs{4mm}

 Let us consider the following  important example. Indeed, we claim to compute the fundamental function
of an ellipsoid relative the norm in anisotropic Lebesgue spaces.\par

 Some additional notations.   $ \vec{1} = (1,1,\ldots,1); $ and for the $  d \ -  $ dimensional vector
$ \vec{p} = p = (p_1, p_2, \ldots, p_{d-2}, p_{d-1},p_d), $
where $ d \ge 2,  $ we define its right-hand side truncation

$$
\vec{p}_{(d)} = p_{(d)} \stackrel{def}{=} (p_1, p_2, \ldots, p_{d-2}, p_{d-1}).
$$

 Let  $  a = \vec{a} = (a_1, a_2, \ldots, a_d) $ be numerical  $  d  \  - $ dimensional vector
with positive entries $ a_i > 0; \hspace{4mm}  d = 1,2,\ldots. $  Define the ellipses (ellipsoids)

$$
E_a = E_{\vec{a}} = \left\{ x = (x_1,x_2,\ldots,x_d): \sum_{i=1}^d \frac{x_i^2}{a_i^2} \le 1  \right\}, \eqno(5.6)
$$

$$
E_a(R) = E_{\vec{a}}(R) = \left\{ x = (x_1,x_2,\ldots,x_d): \sum_{i=1}^d \frac{x_i^2}{a_i^2} \le R^2 \right\} =
E_{R \cdot \vec{a}}, \eqno(5.6a)
$$
so that the ordinary Euclidean centered unit $  d \ - $ dimensional ball $ B $ is equal to the ellipsoid $  B = E_{1,1,\ldots,1} $ and
the ordinary Euclidean centered ball $ B(R)  $ with radii $  R  $ is equal to the ellipsoid  $  B(R) = E_{1,1,\ldots,1}(R).  $\par

\vs{3mm}

 Denote also for simplicity

$$
\theta(\vec{p}, \vec{a}) = \theta^{(d)}(\vec{p}, \vec{a})  = \phi_{\vec{p}} \left(E_{\vec{a}} \right),  \
\theta(\vec{p}) = \theta^{(d)}(\vec{p}) \stackrel{def}{=} \theta(\vec{p}, \vec{1}); \eqno(5.7)
$$
then obviously

$$
\theta^{(d)}(\vec{p}, \vec{a})= \theta^{(d)}(\vec{p}) \cdot \prod_{i=1}^d a_i^{1/p_i} = \theta(\vec{p}, \vec{1}) \prod_{i=1}^d a_i^{1/p_i} \eqno(5.8)
$$
and

$$
\theta(\vec{p}, \vec{a};R) :=
\phi_{\vec{p}} \left(E_{\vec{a}(R)} \right) = \theta^{(d)}(\vec{p}) \cdot \prod_{i=1}^d a_i^{1/p_i} \cdot R^{\sum_i 1/p_i}. \eqno(5.8a)
$$

 We derive after some computations $  \theta^{(1)}(p_1) = 2^{1/p_1}, $

$$
\theta^{(2)}(p_1,p_2) =  \theta^{(1)}(p_1) \ B^{1/p_2} (1/2, 1 + p_2/(2p_1)),
$$
where $  B(\alpha,\beta) $ denotes usually beta function.\par
 Moreover, we can deduce the following recurrent relation

$$
\frac{\theta^{(d+1)}(\vec{p})}{\theta^{(d)}( \vec{p}_{(d+1)} ) } = Z_{d+1}(p_1,p_2, \ldots, p_{d}, p_{d+1}) = Z_{d+1},
$$
where $ \vec{p} \in R^{d+1},  \  Z_1(p_1) = 2^{1/p_1}, $

$$
Z_{d+1}(p_1,p_2, \ldots, p_{d}, p_{d+1}) \stackrel{def}{=}
B^{1/p_{d+1}} \left( \frac{1}{2}, 1 + \frac{p_{d+1}}{2} \ \left( \frac{1}{p_1} +  \frac{1}{p_2} + \ldots + \frac{1}{p_{d}}  \right)  \right). \eqno(5.9)
$$
Therefore

$$
\theta^{(d)}(p_1,p_2,\ldots,p_d) = Z_1(p_1)\cdot Z_2(p_1,p_2) \cdot \ldots Z_d(p_1,p_2,\ldots, p_d) =
$$

$$
\stackrel{def}{=} W_d(p_1,p_2,\ldots, p_d) = W_d(\vec{p}). \eqno(5.10)
$$
 For instance,

$$
\theta^{(3)}(p,q,r) = 2^{1/p} \ B^{1/q}(1/2, 1 + q/(2p)) \ B^{1/r}(1/2, 1 + (r/2)(1/p + 1/q)),
$$

$$
\theta^{(4)}(p,q,r,s) = \theta^{(3)}(p,q,r) \cdot B^{1/s}(1/2, 1 + (s/2)(1/p + 1/q + 1/r))
$$
etc.\par

\vspace{3mm}

Correspondingly

$$
\theta(\vec{p}, \vec{a};R) :=
\phi_{\vec{p}} \left(E_{\vec{a}(R)} \right) = W_{(d)}(\vec{p}) \cdot \prod_{i=1}^d a_i^{1/p_i} \cdot R^{\sum_i 1/p_i}. \eqno(5.11)
$$

\vspace{3mm}

{\bf Remark 5.4.} It is easily to verify that when $  \vec{a} = \vec{p} = \vec{1},  $ the expression (5.11) gives the Euclidean
volume of $ d \ -  $ dimensional ball with radii $ R. $

\vspace{3mm}

{\bf Remark 5.5.} Obviously,  the expression for $ \theta(\vec{p}, \vec{a};R) $ does not dependent from the  center of our ellipsoid.  Namely,
the fundamental function of the non-centered ellipsoid of the form

$$
\tilde{E}_a = \left\{ x = (x_1,x_2,\ldots,x_d): \sum_{i=1}^d \frac{(x_i - x_i^o)^2}{a_i^2} \le R^2  \right\}
$$
is at the same as in the case when $ x_i^o = 0. $\par

\vspace{3mm}

{\bf Remark 5.6.} It is no hard to calculate the fundamental function for parallelepiped

$$
Q = Q(x_1^o, x_2^o, \ldots, x_d^o;  x_1^0+ \delta_1, x_2^o+\delta_2, \ldots, x_d^o + \delta_d) =
$$

$$
\{ \vec{x}: x_j^0 \le x_j \le x_j^o + \delta_j \}, \ \delta_j \in (0,\infty).
$$
 We conclude using  the fact that the indicator function $ I_Q(x) $ is factorable

$$
\phi_{\vec{p}}(Q) = \prod_{j=1}^d \delta_j^{1/p_j}.
$$

\vspace{4mm}

{\bf A. Anisotropic  Lebesgue-Riesz spaces.}
\vspace{3mm}

 Let us return to the announced problem of calculation of the norm of multivariate dilation operator in the mixed
(anisotropic) Lebesgue spaces. Namely,

$$
 R^d  = \otimes_{j=1}^l R^{m_j}, \ \sum_{j=1}^l m_j = d,
$$

 $$
 p = \vec{p} = (p_1, p_2, . . . , p_l), \ p_l \ge 1.
 $$

 Let also $ \vec{A} = A = (A_1, A_2, \ldots, A_l), $ where $ A_j: R^{m_j} \to R^{m_j} $ are linear non - degenerate
 $ \det(A_j) \ne 0 $ operators (matrices), so that $  \vec{A} $ is matrix tensor. \par

Let also $  f(\vec{x}) = f(x_1,x_2,  \ldots, x_l), \ x_j \in R^{m_j} $   be measurable function from the anisotropic space
$ L_{\vec{p}}.  $  We consider the following tensor  dilation operator

$$
V_{\vec{A}}[f] \stackrel{def}{=} f(A_1 x_1, A_2 x_2, \ldots, A_l x_l).  \eqno(5.12)
$$
 Denote

$$
\Lambda_{\vec{p}}(\vec{A}) = \prod_{j=1}^l |\det{A_j}|^{- 1/p_j}. \eqno(5.13)
$$

\vspace{4mm}

{\bf Proposition 5.1.}

$$
|V_{\vec{A}}[f] |_{\vec{p}} \le \Lambda_{\vec{p}}(\vec{A}) \ |f|_{\vec{p}}, \eqno(5.14)
$$
where the equality is attained if for example, the function $  f(\cdot) $ is factorable:

$$
f(x) = \prod_{j=1}^l g_j(x_j), \ g_j(\cdot) \in L_{p_j}(R^{m_j}). \eqno(5.15)
$$

\vspace{4mm}

{\bf B. Anisotropic Grand Lebesgue-Riesz spaces.}

\vspace{3mm}

 Let $ Q $ be convex (bounded or not) subset of the set $ \otimes_{j=1}^l [1,\infty]. $
 Let $ \psi = \psi(\vec{p}) $ be continuous in an interior $ Q^0 $ of the set $ Q $
strictly  positive  function such that

$$
\inf_{\vec{p} \in Q^0}  \psi(\vec{p}) > 0; \ \inf_{\vec{p} \notin Q^0}  \psi(\vec{p}) = \infty.
$$

 We denote the set all of such a functions by $ \Psi(Q). $ \par
The  Anisotropic Grand Lebesgue Spaces $ AGLS = AGLS(\psi) $
 space consists by definition on all the measurable functions

$$
f:  \otimes_{j=1}^l R^{m_j} ( =  R^d) \to R
$$
with finite (mixed) norm

$$
||f||AGLS\psi = \sup_{\vec{p} \in Q^0} \left[ \frac{|f|_{\vec{p}}}{\psi(\vec{p} )} \right]. \eqno(5.16)
$$
  These spaces appear (and investigated) at first (perhaps) in the articles \cite{Ostrovsky101}, \cite{Ostrovsky110}; therein are described
also some possible its applications.  \par

\vs{3mm}

 As before, the (multivariate) fundamental function $  \phi_{AGLS\psi}(D), \ D \subset R^d $ in these space
 can be defined as follows:

$$
\phi_{AGLS\psi}(D) \stackrel{def}{=} ||I_D||AGLS\psi. \eqno(5.17)
$$

\vs{3mm}

 Assume again that the set $  D  $ is direct (Cartesian) product of the (measurable) sets $  F_j: $

$$
D = \otimes_{j=1}^g F_j, \ F_j \subset R^{m_j}.
$$

 Assume in addition  that the function $ \psi(\vec{p})  $ is factorable:

$$
\psi(\vec{p}) = \prod_{j=1}^l \psi_j(\vec{p_j})
$$
and that  the domain $  G  $ is also factorable:

$$
G = \otimes_{j=1}^g G_j, \ G_j \subset R_+^{m_j}.
$$

 Since  the indicator function $ I_D $ is also factorable, we deduce by means of the formula of factorization

$$
\phi_{AGLS\psi}(D) = \prod_{j=1}^g \phi_{ AGLS\psi_j}(F_j).  \eqno(5.18)
$$

\vs{3mm}

 In order to formulate (and prove) the main result of our report, we need to introduce some new preliminary notations.

$$
\delta_j := |\det A_j|^{-1/m_j}, \ j = 1,2, \ldots,l;
$$

$$
K_j^{m_j} := \otimes_{s = 1}^{m_j} [0,\delta_j], \hspace{5mm} K = K(\vec{m}, \vec{p}) := \otimes_{j=1}^l K_j^{m_j}.
$$
so that the set $ K_j^{m_j} $ is a cube of the volume (measure) $ |\det A_j|^{-1} $  in the correspondent space $ R^{m_j}. $\par
 Let further the functions $ \psi = \psi(\vec{p}), \ \zeta = \zeta(\vec{p}), \ \vec{p} \in D $
 be two functions from  certain non - trivial domain $  \Psi(D). $
Define a new functions from this class $  \Psi(D) $

$$
\nu(\vec{p}) = \psi(\vec{p}) \cdot \zeta(\vec{p}).
$$

\vs{4mm}

{\bf Theorem 5.1.}

$$
||V_{\vec{A}}[f]||G\nu  \le \phi(G\zeta, K(\vec{m}, \vec{p}) ) \ ||f||G\psi, \eqno(5.19)
$$
and the last estimate is also in  general case non-improvable, for example  when  the matrix tensor $  \vec{A}  $ is diagonal and
the function $ \psi(\cdot)  $ is natural for the {\it factorable} function $ f: \ \psi(\vec{p}) = |f|_{\vec{p}} < \infty, \ p \in D.  $ \par

\vs{3mm}

{\bf Proof.} Of course, we can and will suppose without loss of generality $  f \in G\psi;  $ in opposite case it is nothing to prove.
Moreover,  it is reasonable to assume $ ||f||G\psi = 1.  $ Then

$$
|f|_{\vec{p}} \le \psi(\vec{p}), \ p \in D.
$$

 We use the inequality (5.14):

$$
|V_{\vec{A}}[f]|_{\vec{p}} \le  \Lambda_{\vec{p}}(\vec{A}) \ |f|_{\vec{p}} \le \Lambda_{\vec{p}}(\vec{A}) \
  \psi(\vec{p}), \ p \in D,
$$
or equally

$$
\frac{|V_{\vec{A}} [f]|_{\vec{p}}}{\nu(\vec{p})}
 \le \frac{\Lambda(\vec{A})}{\zeta(p)}. \eqno(5.20)
$$

 Let us take the supremum over $ p \in D $ from both the sides of the last inequality, taking into account
the direct definitions of norm and fundamental function for Anisotropic Grand Lebesgue Spaces:

$$
||V_A[f]||G\nu \le \phi(G\zeta, K(\vec{m}, \vec{p}) )=
\phi(G\zeta, K(\vec{m}, \vec{p}) ) \ ||f||G\psi, \eqno(5.21)
$$
Q.E.D.

\vspace{4mm}

 \section{ Concluding remark. Examples and counterexamples.   }

\vspace{4mm}

{\bf A. An example.}\\

\vs{3mm}

 Let ua show that the condition $  \det(A) \ne 0 $ in the theorem 3.1 (and in another ones) is essential. Let $ f = f(x,y), (x,y) \in R^2 $ be
measurable non-negative factorable  function

$$
f(x,y) = g(x) \cdot g(y),
$$
 where $ g(\cdot) \in L_p(R), \ \lim_{y \to 0} g(y) = \infty, $
 and consider the linear degenerate operator $ A  \ $ (matrix $  2 \times 2) $ such that $  A(x,y) = (x,0), $ so that

$$
f(A\cdot(x,y)) := f(x,0).
$$
  So, the operator $  A  $ is the coordinate projections.\par

  Evidently, $ f(\cdot, \cdot) \in L_p(R^2),  $  but
  $$
  f(A(x, y)) = f(x,0) \notin L_p(R^2).
  $$

  Note that the coordinate projections on matrix weighted $ L_p \ - $ spaces,  for instance, Hilbert's transform
 is investigated in the recent article \cite{Nielsen1}.\par

\vspace{3mm}

{\bf B. Particular case: Exponential Orlicz spaces.}\\

\vspace{3mm}

 It is known, see \cite{Kozatchenko1}, \cite{Ostrovsky1}, \cite{Liflyand1}, \cite{Ostrovsky7}, that the so-called
{\it exponential} Orlicz's spaces, for example, the Orlicz's spaces with correspondent Young function

$$
\Phi(u) = \Phi^{(\lambda)}(u) =
e^{C |u|^{\lambda} } - 1, \ \lambda, C = \const > 0
$$
are particular case of Grand Lebesgue Space with correspondent $ \psi = \psi_{\Phi} \ - $ function, for example,

$$
\psi_{\Phi^{(\lambda)}}(p) := \psi^{(\lambda)}(p) = p^{1/\lambda}, \ p \ge 1.
$$

 Therefore, the theorem 3.1 may be applied to these Orlicz spaces. \par
 See also \cite{Estaremi1}, \cite{Estaremi2}, \cite{Gupta1} etc. \par
The correspondent fundamental function $  \phi(G\psi^{(\lambda)}, \delta) $ is  calculated and estimated in
\cite{Ostrovsky2}.\par

\vspace{3mm}

{\bf C. Periodical case.}\\

\vspace{3mm}

  At the same results may be derived in the case when $  X = (-\pi, \pi)^d $ (case of torus), where the algebraic operations
 are understood  $  \mod(2 \pi); $  or more generally when

$$
 X = (-\pi, \pi)^{d_1} \otimes R^{d_2}, \ d_1, d_2 = 1,2,\ldots.
$$

\vspace{3mm}

{\bf D. Example to the third section.}\\

\vspace{3mm}

 Let us consider the $ G\psi \ - $ space $ \tilde{G} =  \tilde{G}(a,h; \alpha,\beta)  $ over real
line $  R^1 $ with the following $  \psi \ - $ function

$$
\tilde{\psi} = \psi(a; \alpha,\beta; p) = (p-a)^{-\alpha}, \ p \in (a,h);
$$

$$
\tilde{\psi}(p) = \psi(a; \alpha,\beta; p) = p^{\beta}, \ p \in (h, \infty),
$$
where

$$
\alpha, \beta = \const > 0, \ a = \const \ge 1,
$$
the value $ h = h(a,\alpha, \beta) > a $ is the unique positive solution of an equation

$$
(h- a)^{-\alpha} = h^{\beta},
$$
so that the function $ p \to \tilde{\psi}(p), \ p \in (a,\infty) $ is continuous.\par
 This  space  does not coincides  with  the known rearrangement invariant spaces: Lorentz, Orlicz,
Marcinkiewicz etc., see \cite{Kozatchenko1}, \cite{Ostrovsky1}, \cite{Ostrovsky2}, \cite{Liflyand1}, \cite{Ostrovsky7}. \par

The correspondent fundamental function $ \tilde{\phi}(\delta) = \tilde{\phi}((a;\alpha,\beta,\delta)$ obeys a following
asymptotical behavior:

$$
\tilde{\phi}(\delta) \sim \beta^{\beta} \ |\ln \delta|^{-\beta}, \ \delta \to 0+;
$$

$$
\tilde{\phi}(\delta) \sim  \left(\frac{a^2 \alpha}{e} \right)^{\alpha} \ \delta^{1/\alpha} \ (\ln \delta)^{-\alpha}, \ \delta \to \infty,
$$
 see \cite{Ostrovsky2}. \par
It remains to apply the inequality (3.3). \par

\vspace{3mm}

{\bf E. Possible generalisations.}\\

\vspace{3mm}

 It may be investigated analogously to the $  5^{th} $ section the multivariate  weight case as well as the affine linear
non-centered transform  of the form

$$
 f \to f(A\cdot x + b).
$$

\vspace{3mm}

{\bf F. About compactness of dilation operator.}\\

\vspace{3mm}

 Of course, in general case the dilation operator $  V_A[\cdot] $ acting  from one GLS to  suitable ones, is
non- compact, even in the very simple case $  A = I \ - $  unit operator: $  V_A[f] = f. $\par
 Let us consider an opposite case. We borrow the notations, conditions and proposition of theorem 3.1. \par
Let also $ \psi_1(\cdot), \ \psi_2(\cdot) $  be two function from the set $ \Psi(a,b), \ 1 \le a < b \le \infty. $
We recall the following relation definition:

$$
\psi_1 << \psi_2 \ \Leftrightarrow  \ \lim_{\psi_2(p) \to \infty} \frac{\psi_1(p)}{\psi_2(p)} = 0,
$$
see \cite{Liflyand1},  \cite{Ostrovsky2}. It is proved in these articles that in this case
the GLS space $  G\psi_1 $ is compact embedded in the space  $  G\psi_2. $ \par

 We deduce as a slight corollary: let the new function $ \theta = \theta(p) $ from this set  $ \Psi(a,b) $
be such that $ \theta(\cdot) << \nu(\cdot).  $  Then we derive under conditions of theorem 3.1
that the operator $ V_A[\cdot] $ acting from the space $ G\psi $ into the space $ G\theta $
is compact. \par

\end{document}